\documentclass[12pt]{article}

\usepackage{amsmath,amsfonts,amssymb,lmodern,geometry,enumerate}
\usepackage[T1]{fontenc}
\usepackage[english]{babel}
\usepackage{dsfont}
\usepackage{amsthm,stmaryrd}
\usepackage{color}
\usepackage{bm,bbm}
\usepackage{mathrsfs,url,color}
\usepackage{textcase}
\usepackage{graphicx}
\usepackage{wrapfig}

\newtheorem{thm}{Theorem}[section]

\newtheorem{defi}[thm]{Definition}
\newtheorem{lma}[thm]{Lemma}

\newtheorem{re}{Remark}[section]
\newtheorem{ap}{Application}

\newcommand{\Pro}{\mathbb{P}} 
\newcommand{\prob}{\Pro}
\newcommand{\E}{\mathbb{E}}
\newcommand{\mean}{\E}
\newcommand{\real}{\mathbb{R}}

\newcommand{\Z}{\mathbb{Z}}

\newcommand{\bone}{{\bf 1}}

\def\Ceka{\v Cekan\-avi\v cius}

\def\equald{\stackrel{\mbox{\scriptsize{d}}}{=}}

\def\im{{\mbox{\sl i}}}
\def\ime{{\mbox{\sl\scriptsize i}}}

\def\scrB{{\mathscr{B}}}

\def\var{{\rm Var}}

\def\ignore#1{}
\def\Ref#1{(\ref{#1})}
\def\qed{\hfill\hbox{${\vcenter{\vbox{
        \hrule height 0.4pt\hbox{\vrule width 0.4pt height 6pt
        \kern5pt\vrule width 0.4pt}\hrule height 0.4pt}}}$}}

\newcounter{con}
\stepcounter{con}
\newcommand{\qcon}[1]{\addtocounter{con}{1}}

\numberwithin{equation}{section}

\ignore{\makeatletter
\renewcommand\section{\@startsection {section}{1}{\z@}%
{-3.5ex \@plus -1ex \@minus -.2ex}%
{1.3ex \@plus.2ex}%
{\center\small\sc\MakeTextUppercase}}
\def\subsection#1{\@startsection {subsection}{2}{0pt}%
{-3.5ex \@plus -1ex \@minus -.2ex}%
{1ex \@plus.2ex}%
{\bf\mathversion{bold}}{#1}}
\def\subsubsection#1{\@startsection{subsubsection}{3}{0pt}%
{\medskipamount}%
{-10pt}%
{\normalsize\itshape}{\kern-2.2ex. #1.}}
\makeatother}

\begin{document}

\title{\sc\bf\large\MakeUppercase{
A dichotomy for CLT in total variation 
}}

\author{Aihua Xia\thanks{School of Mathematics and Statistics, the University of Melbourne, 
Parkville, VIC 3010, Australia; aihuaxia@unimelb.edu.au;
work supported by Australian Research Council Discovery Grant DP150101459.} \\ \it University of Melbourne}

\date{\today}

\maketitle
\vskip-1cm
\begin{abstract}
Let $\eta_i$, $i\ge 1$, be a sequence of independent and identically distributed random variables with finite third moment,
and let $\Delta_n$ be the total variation distance between the distribution of $S_n:=\sum_{i=1}^n\eta_i$ and the normal distribution with the same mean and variance. In this note, we show the dichotomy that either $\Delta_n=1$ for all $n$ or $\Delta_n=O\left(n^{-1/2}\right)$.
\end{abstract}



\section{Introduction and the main result}
The Berry-Esseen Theorem (Berry~\cite{Berry41} and Esseen~\cite{Esseen42}) states that if $\eta_i$, $1\le i\le n,$ are independent and identically distributed (iid) random variables with mean 0 and variance 1,
$Y_n=\frac{\sum_{i=1}^n\eta_i}{\sqrt{n}}$, $Z\sim N(0,1)$, where $\sim$ denotes ``is distributed as'', then 
$$d_K(Y_n,Z)\le \frac{C\mean|\eta_1|^3}{\sqrt{n}},$$
where $d_K$ is the Komogorov distance: for two random variables $X_1$ and $X_2$ with distributions $F_1$ and $F_2$,
$$d_K(X_1,X_2):=d_K(F_1,F_2):=\sup_{x\in\real}|F_1(x)-F_2(x)|.$$
The Kolmogorov distance $d_K(F_1,F_2)$ measures the difference between the distribution functions $F_1$ and $F_2$, but it does not
tell much about the difference between the probabilities $\prob(X_1\in A)$ and $\prob(X_2\in A)$ for a non-interval Borel set $A\subset\real$, e.g., $A=\cup_{i\in \Z}(2i-0.1,2i+0.1)$, where $\Z$ denotes the set of all integers. Such difference is reflected in the total variation distance $d_{TV}(F_1,F_2)$ defined by
$$d_{TV}(X_1,X_2):=d_{TV}(F_1,F_2):=\sup_{A\in\scrB(\real)}|F_1(A)-F_2(A)|,$$
where $\scrB(\real)$ denotes the Borel $\sigma$-algebra on $\real$ and $F_i(A):=\int_AdF_i(x)$. The definition is equivalent to 
\begin{equation}d_{TV}(F_1,F_2)=\frac12\sup_f\left|\int f(x)dF_1(x)-\int f(x) dF_2(x)\right|,\label{totalvareq}\end{equation}
where the supremum is taken over all measurable functions $f$ on $(\real,\scrB(\real))$ bounded by 1.

Although central limit theorems in the total variation have been studied in some special circumstances (see, e.g., \cite{DF87,MM07}), it is generally believed that the total variation distance is too strong for normal approximation (see, e.g., \Ceka~\cite{Ceka00}, Chen and Leong~\cite{CL10}, Fang~\cite{Fang14}). For example, the total variation distance between any binomial distribution and any normal distribution is always 1. To recover central limit theorems in the total variation, a common approach is to discretize the distribution of interest and approximate it with a simple discrete distribution, e.g.,
 translated Poisson (R\"ollin~\cite{Rollin05,Rollin07}), centered binomial (R\"ollin~\cite{Rollin08}), discretized normal (Chen and Leong~\cite{CL10}, Fang~\cite{Fang14}) and a family of polynomial type distributions (Goldstein and Xia~\cite{GX06}). The multivariate versions of these approximations are investigated by Barbour, Luczak and Xia~\cite{BLX15}.

By discretizing a distribution $F$ of interest, we essentially group the probability of an area and put it at one point in the area, hence the information of $F(A)$ for a general set $A\in\scrB(\real)$ is completely lost. 
In this note, we consider the normal approximation in the total variation to the sum of iid random variables with finite second moment.

The Lebesgue decomposition theorem \cite[p.~134]{Halmos74} ensures that any distribution function $F$ on $\real$ can be represented as
$$F=(1-\alpha_F) F_s+\alpha_FF_a,$$
where $\alpha_F\in[0,1]$, $F_s$ and $F_a$ are two distribution functions such that, with respect to the Lebesgue measure on $\real$, $F_a$ is absolutely continuous and $F_s$ is singular.

\begin{defi} A distribution function $F$ on $\real$ is said to be non-singular if $\alpha_F>0$.
\end{defi}

In other words, $F$ is non-singular if and only if there exists a sub-probability measure $F_0\ne 0$ with a density $f_0$ such that
\begin{equation}F_0(A)=\int_A f_0(x)dx,\mbox{ for all }A\in\scrB(\real).\label{defi1} \end{equation}

\ignore{The family of non-singular distributions is huge, including uniform, beta, exponential, triangular, extremal, spherically symmetric~\cite{Boisbunon12,DF87} and zero adjusted distributions~\cite[Ch.~5]{RSH14}.}

\begin{thm}\label{thm1} Let $\eta_i$, $i\ge 1,$ be iid with finite second moment, define $S_n:=\sum_{i=1}^n\eta_i$ and 
$\Delta_n=d_{TV}(S_n,Z_n)$, where $Z_n\sim N(\mean S_n,\var(S_n))$.
 The following are equivalent:
 \begin{description}
\item{(i)} There exists a finite integer $n_0$ such that $\Delta_{n_0}<1$.
\item{(ii)} There exists a finite integer $n_0$ such that $S_{n_0}$ is non-singular.
\item{(iii)} $\Delta_n=o\left(1\right)$.
\end{description}
Furthermore, if $\eta_1$ has finite third moment, then (i)-(iii) are also equivalent to
 \begin{description}
\item{(iv)} $\Delta_n=O\left(n^{-1/2}\right)$.
\end{description}
\end{thm}

\begin{re} Theorem~\ref{thm1} says that, for the iid sequence $\{\eta_i:\ i\ge 1\}$ with finite third moment, we have the dichotomy that either $S_n$ is singular so that $\Delta_n=1$ for all $n$ or it can be approximated by the normal distribution with the same mean and variance in the total variation with convergence speed no less than $O\left(n^{-1/2}\right)$. Curiously, this phenomenon may be related to the Kolmogorov's zero-one law.
\end{re}

It is possible to generalize some parts of Theorem~\ref{thm1} to non-identically distributed random variables but the formulation of such generalizations is typically complicated. By focusing on the most important case, we aim to keep the paper reader-friendly and to deliver a clear and concise message. 

The proof of Theorem~\ref{thm1} is based on Stein's method for normal approximation and the estimate of $d_{TV}(S_n,S_n+\gamma)$. Generally speaking, the easiest metric that Stein's method for normal approximation can handle is the Wasserstein metric. Much more effort is needed to 
achieve an error bound for the Kolmogorov distance enjoying the same order as that for the Wasserstein distance bound \cite{Chen98,CGS11}. The way that Stein's method for normal approximation is used in this paper seems to be unexplored. 
In the context of Poisson and other discrete distribution approximations, this approach is well studied (see, e.g., \cite{Xia97,BX99,BLX15}).

\section{The proof of Theorem~\ref{thm1}}

We start with a few technical lemmas.

\begin{lma} \label{lma1} Assume $\xi_1,\dots,\xi_n$ are iid random variables having the triangular density function
\begin{equation}\kappa_a(x)=\left\{\begin{array}{ll}
\frac1a\left(1-\frac {|x|}a\right),&\mbox{ for }|x|\le a,\\
0,&\mbox{ for } |x|>a,
\end{array}\right.\label{lma1.1}\end{equation}
where $a>0$. Let $T_n=\sum_{i=1}^n\xi_i$. Then for any $\gamma>0$,
\begin{equation}d_{TV}(T_n,T_n+\gamma)\le \frac{\gamma}a\left\{\sqrt{\frac{3}{\pi n}}+\frac{2}{(2n-1)\pi^{2n}}\right\}.\label{lma1.2}\end{equation}
\end{lma}

\noindent{\it Proof.} For convenience, we write $G_n$, $g_n$ and $\psi_n$ as the distribution, density and characteristic functions of $T_n$ respectively. It is well-known that the triangular density $\kappa_a$ has the characteristic function $\psi_1(s)=\frac{2(1-\cos(as))}{(as)^2}$, which gives $\psi_n(s)=\left(\frac{2(1-\cos(as))}{(as)^2}\right)^n$. Using the fact that the convolution of two symmetric unimodal distributions on $\real$ is unimodal \cite{Wintner38}, we can conclude that the distribution of $T_n$ is unimodal and symmetric. This ensures that 
\begin{equation}d_{TV}(T_n,T_n+\gamma)=\sup_{x\in\real}|G_n(x)-G_n(x-\gamma)|=\int_{-\gamma/2}^{\gamma/2}g_n(x)dx.\label{lma1.3}\end{equation}
Applying the inversion formula, we have
\begin{eqnarray*}
g_n(x)&=&\frac1{2\pi}\int_\real e^{-\ime sx}\psi_n(s)ds=\frac1{2\pi}\int_\real \cos(sx)\psi_n(s)ds\nonumber\\
&=&\frac1{a\pi}\int_0^\infty \cos(sx/a)\left(\frac{2(1-\cos s)}{s^2}\right)^nds,
\end{eqnarray*}
where $\im=\sqrt{-1}$ and the second equality is due to the fact that $\sin(sx)\psi_n(s)$ is an odd function. Obviously, $g_n(x)\le g_n(0)$ so we need to establish an upper bound for $g_n(0)$. A direct verification gives
$$0\le \frac{2(1-\cos s)}{s^2}\le e^{-\frac{s^2}{12}}\mbox{ for }0\le s\le 2\pi,$$
which implies
\begin{eqnarray}
g_n(0)&\le&\frac1{a\pi} \left\{\int_0^{2\pi}e^{-\frac{ns^2}{12}}ds+\int_{2\pi}^\infty \left(\frac 4{s^2}\right)^nds\right\}\nonumber\\
&\le& \frac1{a\pi\sqrt{n}}\int_0^\infty e^{-\frac{s^2}{12}}ds+\frac{2}{a(2n-1)\pi^{2n}}\nonumber\\
&=&\frac1a\sqrt{\frac{3}{\pi n}}+\frac{2}{a(2n-1)\pi^{2n}}.\label{lma1.5}
\end{eqnarray}
Now, combining \Ref{lma1.5} with \Ref{lma1.3} gives \Ref{lma1.2}.\qed 

We denote the convolution by $\ast$ and write $F^{k\ast}$ as the $k$-fold convolution of the function $F$ with itself. 

\begin{lma}\label{lma2} If $F$ is a non-singular distribution, then there exist $a>0$, $u\in\real$, $\theta\in(0,1]$ and a distribution function $H_2$ such that 
$$F^{2\ast}=(1-\theta)H_2+\theta H_1\ast\delta_u,$$
where $H_1$ is the distribution of the triangular density $\kappa_a$ in \Ref{lma1.1} and $\delta_u$ is the Dirac measure at $u$.
\end{lma}

\begin{wrapfigure}{r}{0.5\textwidth}
  \vspace{-20pt}
  \begin{center} 
   \includegraphics[trim = 0mm 0mm 0mm 0mm, clip,width=0.5\textwidth]{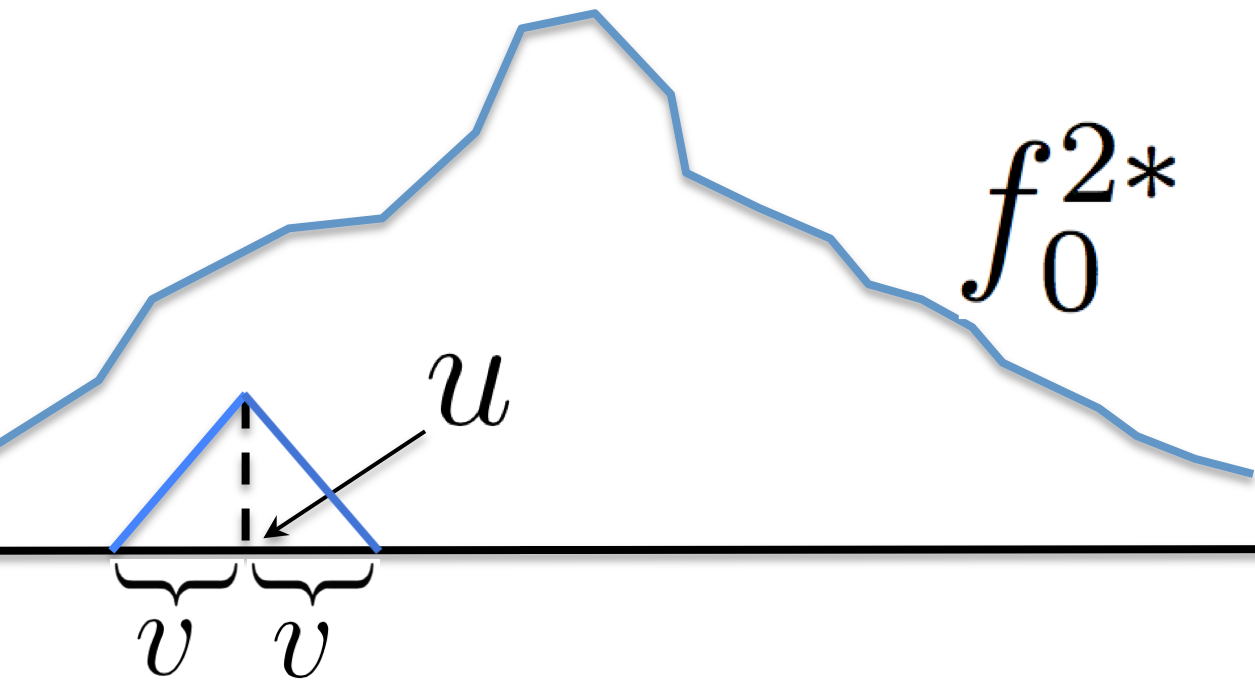}
 \vspace{-20pt}
  \caption{Existence of $u$ and $v$} %
  \label{figureone}
  \end{center}
  \vspace{-10pt}
 \end{wrapfigure}
\noindent{\it Proof.} Since $F$ is non-singular, we can find $f_0$ satisfying \Ref{defi1} and $\int_\real f_0(x)dx>0$. 
Without loss, we assume that $f_0$ is bounded with bounded support. Then $f_0^{2\ast}$ is continuous \cite[p.~79]{Lindvall92}. Referring to Figure~\ref{figureone}, since $f_0^{2\ast}\not\equiv0$, we can find  $u\in\real$ and $v>0$ such that $f_0^{2\ast}(u)>0$ and $\min_{x\in[u-v,u+v]}f_0^{2\ast}(x)\ge \frac12f_0^{2\ast}(u)=:b$. Let $\theta=vb$ and $a=\sqrt{v/b}$, $H_2=\frac1{1-\theta}(F^{2\ast}-\theta H_1\ast\delta_u)$, the claim follows. \qed

Lemma~\ref{lma2} says that $F^{2\ast}$ is the distribution function of $(X_1+u)X_3+X_2(1-X_3)$, where $X_1\sim H_1$, $X_2\sim H_2$, $X_3\sim{\rm Bernoulli}(\theta)$ are independent random variables.

\begin{lma}\label{lma3} With the setup in Theorem~\ref{thm1}, for any $\gamma>0$, we have
\begin{equation}
d_{TV}(S_n,S_n+\gamma)\le (\gamma\vee 1) O(n^{-1/2}),\label{lma3.1}
\end{equation}
where $O(n^{-1/2})$ does not depend on $\gamma$.
\end{lma}

\noindent{\it Proof.} Let $m=\lfloor n/2\rfloor$, the integer part of $n/2$. By Lemma~\ref{lma2}, we can construct independent random variables $X_{ij}$, $1\le i\le 3$, $1\le j\le m$ such that $X_{1j}\sim H_1$, $X_{2j}\sim H_2$, $X_{3j}\sim{\rm Bernoulli}(\theta)$ and $S_{2m}\equald S_m':=\sum_{j=1}^m[(X_{1j}+u)X_{3j}+X_{2j}(1-X_{3j})]$. Clearly $I:=\sum_{j=1}^mX_{3j}$ follows Binomial$(m,\theta)$ and given $I=k$, the conditional distribution of 
$S_m'$ is $(H_1\ast\delta_u)^{k\ast}\ast H_2^{(m-k)\ast}$. In other words, 
\begin{equation}S_{2m}\sim\sum_{k=0}^m\prob(I=k)(H_1+\delta_u)^{k\ast}\ast H_2^{(m-k)\ast}.\label{lma3.2}\end{equation}
Let $A-\gamma:=\{x-\gamma:\ x\in A\}$. Using the fact that for distribution functions $G_i$, 
$d_{TV}(G_1\ast G_3,G_2\ast G_3)\le d_{TV}(G_1,G_2)$ in the last two inequalities below, we obtain
from \Ref{lma3.2} that
\begin{eqnarray}
&&d_{TV}(S_n,S_n+\gamma)\le d_{TV}(S_{2m},S_{2m}+\gamma)\nonumber\\
&&=\sup_{A\in\scrB(\real)}\left|\sum_{k=0}^m\prob(I=k)\left\{(H_1+\delta_u)^{k\ast}\ast H_2^{(m-k)\ast}(A)-(H_1+\delta_u)^{k\ast}\ast H_2^{(m-k)\ast}(A-\gamma)\right\}\right|\nonumber\\
&&\le\sum_{k=0}^m\prob(I=k)\sup_{A\in\scrB(\real)}\left|(H_1+\delta_u)^{k\ast}\ast H_2^{(m-k)\ast}(A)-(H_1+\delta_u)^{k\ast}\ast H_2^{(m-k)\ast}(A-\gamma)\right|\nonumber\\
&&=\sum_{k=0}^m\prob(I=k)d_{TV}\left((H_1+\delta_u)^{k\ast}\ast H_2^{(m-k)\ast},(H_1+\delta_u)^{k\ast}\ast\delta_\gamma\ast H_2^{(m-k)\ast}\right)\nonumber\\
&&\le\sum_{k=0}^m\prob(I=k)d_{TV}\left((H_1+\delta_u)^{k\ast},(H_1+\delta_u)^{k\ast}\ast\delta_\gamma\right)\nonumber\\
&&=\sum_{k=0}^m\prob(I=k)d_{TV}\left(H_1^{k\ast},H_1^{k\ast}\ast\delta_\gamma\right)\nonumber\\
&&\le d_{TV}\left(H_1^{k_0\ast},H_1^{k_0\ast}\ast\delta_\gamma\right)+\prob(I\le k_0-1).\label{lma3.3}
\end{eqnarray}
If we take $k_0=\lfloor 0.5m\theta\rfloor$, it follows from Lemma~\ref{lma1} and \Ref{lma3.3} that
$$d_{TV}(S_n,S_n+\gamma)\le \gamma O(m^{-1/2})+O(m^{-1}).$$
However, $m=\lfloor n/2\rfloor$, the proof is complete. \qed

\noindent{\it Proof of Theorem~\ref{thm1}.} If $\var(\eta_1)=0$, then $S_n$ is singular and $\Delta_n=1$ for all $n$. Now, we assume $\var(\eta_1)>0$. As the total variation distance is invariant in terms of the linear transformation, without loss, we can further assume that $\eta_i$'s have mean 0 and variance 1. Define $W_n=\frac{S_n}{\sqrt{n}}$, then $\mean W_n=0$ and $\var(W_n)=1$.

(i)$\Rightarrow$(ii): If $S_{n_0}$ is singular, then there exists
an $A\in\scrB(\real)$ with $|A|=\int_A dx=0$ and $\prob(S_{n_0}\in A)=1$. This implies that $\Delta_{n_0}=1$, which contradicts (i).

(ii)$\Rightarrow$(iii): The Stein equation for the standard normal distribution (see  \cite[p.~15]{CGS11}) is
\begin{equation}f'(w)-wf(w)=h(w)-Nh,\label{Stein1}\end{equation}
where $Nh:=\mean h(Z)$ for $Z\sim N(0,1)$. The solution of the Stein equation satisfies (see \cite[p.~16]{CGS11})
$$\|f_h'\|:=\sup_w\left|f_h'(w)\right|\le 2\|h(\cdot)-Nh\|.$$ 
Hence, for $h=\bone_A$ with $A\in\scrB(\real)$, the solution $f_h=:f_A$ satisfies 
\begin{equation}
\|f_A'\|\le 2.\label{thm1.3}
\end{equation}
Regrouping $\eta_i$, $1\le i\le n$, into blocks of $n_0$ random variables and setting $\eta_i'=\sum_{j=(i-1)n_0+1}^{in_0}\eta_j$ if necessary, we may assume that $\eta_i$'s are non-singular. For convenience, we define $W_n':=W_n-\eta_1/\sqrt{n}$ and omit the subindex $A$ in $f_A$, then
\begin{eqnarray}
&&\mean [f'(W_n)-W_nf(W_n)]\nonumber\\
&&=\mean f'(W_n)-n\mean [\eta_1/\sqrt{n}(f(W_n'+\eta_1/\sqrt{n})-f(W_n'))]\nonumber\\
&&=\mean f'(W_n)-\mean \eta_1^2\int_0^1f'(W_n'+u\eta_1/\sqrt{n})du\nonumber\\
&&=\mean [f'(W_n)-f'(W_n')]-\mean \eta_1^2\int_0^1(f'(W_n'+u\eta_1/\sqrt{n})-f'(W_n'))du\nonumber\\
&&=\mean\left\{\mean \left[f'(W_n'+\eta_1/\sqrt{n})-f'(W_n')\big|\eta_1\right]\right\}\nonumber\\
&&\ \ \ -\int_0^1\mean\left\{\mean\left[f'(W_n'+u\eta_1/\sqrt{n})-f'(W_n')\big|\eta_1\right]\eta_1^2\right\}du,\label{thm1.4}
\end{eqnarray}
where, since $\eta_1$ is independent of $W_n'$, the first equality is guaranteed by $\mean \eta_1f(W_n')=\mean \eta_1\mean f(W_n')=0$, and the third equality follows from $\mean \eta_1^2f'(W_n')=\mean \eta_1^2\mean f'(W_n')=\mean f'(W_n')$.  Now, with 
$$d_{n,v}:=d_{TV}(W_n'+|v|/\sqrt{n},W_n')=d_{TV}(S_{n-1}+|v|,S_{n-1}),$$ 
we have from \Ref{totalvareq} that
 \begin{equation}
 \left|\mean \left[f'(W_n'+\eta_1/\sqrt{n})-f'(W_n')\big|\eta_1=v\right]\right|
 \le2\|f'\|d_{n,v}\label{thm1.5}
\end{equation}
and
\begin{equation}\left|\mean\left[f'(W_n'+u\eta_1/\sqrt{n})-f'(W_n')\big|\eta_1=v\right]\right|\le 2\|f'\|d_{n,vu}.\label{thm1.6}
\end{equation}
 We denote the distribution function of $\eta_1$ by $F_\eta$. For $A\in\scrB(\real)$, using \Ref{Stein1} and \Ref{thm1.3}, we obtain from \Ref{thm1.4}, \Ref{thm1.5} and \Ref{thm1.6} that
\begin{equation}
|\prob(W_n\in A)-\prob(Z\in A)|
\le 4\int_\real \left\{d_{n,v}+\left(\int_0^1 d_{n,vu}du\right)v^2\right\}dF_\eta(v).\label{thm1.8}
\end{equation}
Taking supremum over all $A\in\scrB(\real)$ in \Ref{thm1.8}, we get
\begin{equation}\Delta_n\le 4\int_\real \left\{d_{n,v}+\left(\int_0^1 d_{n,vu}du\right)v^2\right\}dF_\eta(v).\label{thm1.9}
\end{equation}
However, for all $s\in\real$, \Ref{lma3.1} ensures that $d_{n,s}\to 0$ as $n\to\infty$. Since $0\le d_{n,s}\le 1$, (iii) follows from the dominated convergence theorem.

(iii)$\Rightarrow$(i) is obvious.

(ii)$\Rightarrow$(iv): It follows from \Ref{thm1.9} and \Ref{lma3.1} that
\begin{eqnarray*}
\Delta_n&\le& O\left(n^{-1/2}\right)\mean (|\eta_n|\vee 1)+O\left(n^{-1/2}\right)\int_0^1\left\{\int_\real (|vu|\vee 1)v^2dF_\eta(v)\right\}du\nonumber\\
&\le& O\left(n^{-1/2}\right)\mean (|\eta_n|\vee 1)+O\left(n^{-1/2}\right)\mean (|\eta_1|\vee 1)^3\nonumber\\
&=&O\left(n^{-1/2}\right),
\end{eqnarray*}
since $\mean|\eta_1|^3<\infty$, concluding the proof. 

(iv)$\Rightarrow$(iii) is also obvious. \qed


\def\ac{{Academic Press}~}
\def\aap{{Adv. Appl. Prob.}~}
\def\ap{{Ann. Probab.}~}
\def\anap{{Ann. Appl. Probab.}~}
\def\eljp{{\it Electron.\ J.~Probab.\/}~} 
\def\jap{{J. Appl. Probab.}~}
\def\jws{{John Wiley $\&$ Sons}~}
\def\ny{{New York}~}
\def\ptrf{{Probab. Theory Related Fields}~}
\def\sp{{Springer}~}
\def\spa{{Stochastic Processes and their Applications}~}
\def\sv{{Springer-Verlag}~}
\def\tpa{{Theory Probab. Appl.}~}
\def\zw{{Z. Wahrsch. Verw. Gebiete}~}

\end{document}